\documentclass[a4paper, 11pt]{article}

\usepackage{amsmath}
\usepackage{amssymb}
\usepackage{graphics}
\usepackage[dvips]{graphicx}
\usepackage{epsfig}
\usepackage[english]{babel}
\usepackage{lineno}
\usepackage{color}
\usepackage{caption}
\usepackage{subcaption}
\usepackage{epstopdf}
\usepackage{tikz}
\usepackage{hyperref}
\graphicspath{ {./figures/} }
\newtheorem{theo}{Theorem}

\newtheorem{prop}{Proposition}

\begin{document}
%\linenumbers*[1]

\title{SIS epidemic propagation on hypergraphs}

\author{\'Agnes Bod\'o$^{1,2}$, Gyula Y. Katona$^{2,3}$, P\'eter L. Simon$^{1,2,\ast}$}

\maketitle

\begin{center}
$^1$ Institute of Mathematics, E\"otv\"os Lor\'and University Budapest, Hungary \\
$^2$ Numerical Analysis and Large Networks Research Group,\\ Hungarian Academy of Sciences, Hungary\\
$^3$ Department of Computer Science and Information Theory, \\ Budapest University of Technology and Economics, Budapest, Hungary
\end{center}

\vspace{1cm}

\begin{abstract}
Mathematical modeling of epidemic propagation on networks is extended to hypergraphs in order to account for both the community structure and the nonlinear dependence of the infection pressure on the number of infected neighbours. The exact master equations of the propagation process are derived for an arbitrary hypergraph given by its incidence matrix. Based on these, moment closure approximation and mean-field models are introduced and compared to individual-based stochastic simulations. The simulation algorithm, developed for networks, is extended to hypergraphs. The effects of hypergraph structure and the model parameters are investigated via individual-based simulation results.
\end{abstract}


\vspace{1cm}


\vspace{1cm}
\begin{flushleft}
$\ast$ corresponding author\\
email: simonp@cs.elte.hu\\
\end{flushleft}

\newpage

%%%%%%%%%%%%%%%%%%%%%%%%%%%%%%%%%%%%%%%%%%%%%%%%%%%%%%%%%%%%%%%%%%%%%%%%%%%%%%%%%%%%%%%%%%
\section{Introduction}
%%%%%%%%%%%%%%%%%%%%%%%%%%%%%%%%%%%%%%%%%%%%%%%%%%%%%%%%%%%%%%%%%%%%%%%%%%%%%%%%%%%%%%%%%%

Spreading processes on networks has a well-developed theory, based on
the simple idea that the higher connectivity of the network enhance
the spreading process. In the language of epidemic propagation this
means that infection pressure on a susceptible individual increases
with the number of infectious neighbours. Typically it is assumed that
the probability of infection is proportional to the number of
infectious neighbours. This idea leaded to different propagation
models from exact master equations \cite{simon2010JMB, van2009} to
different types of mean-field models, such as homogeneous and
heterogeneous pairwise models \cite{Danon, Gleeson, house2011},
effective degree models \cite{Lindquist}, edge based compartmental
models \cite{miller2011edge} and individual based models
\cite{Sharkey08, van2011n}, mentioning only the most widely used
ones. It is also well-known that the community structure has strong
impact on the spread of the epidemic. There are many publications
studying household structure, workplaces, schools.

A real modeling framework of epidemic propagation has to take into
account that (a) the community is built up from small units, such as
households and workplaces and (b) the infection pressure on a
susceptible individual in a unit is not proportional to the number of
infected individuals. For example, the increase of the number of
infected individuals from 5 to 10 in a workplace with 20 individuals
does not necessarily mean that the infection probability is doubled in
that place.  The aim of this paper is to develop the theory of
epidemic propagation on hypergraphs that enables us to model both the
nonlinear dependence of the infection pressure and the community
structure. The modeling paradigm is that the population consists of
nodes of a hypergraph that is given by the hyperedges, which are
simply subsets of the vertex set. For example, each household and each
workplace can be considered as a hyperedge. The usual graphs can be
considered %%%
as the special case of hypergraphs, when each hyperedge
consists of two vertices that form an edge of the graph. In the
widely-used propagation models the rate of infection for a susceptible
node in a household is $n\tau$, where $n$ is the number of its
infectious neighbours in the household and $\tau$ is the per contact
infection rate. In our new model, to be developed in this paper, it is
$f(n)\tau$ with some possibly non-linear function $f$ defined
later. This function describes that the infection pressure is
discounted as the number of infectious individuals is increased. If an
individual belongs to two or more hyperedges, e.g. to a household and
to a workplace, then the rate of infection is the sum of the rates in
each hyperedge, e.g. $(f(n_1)+f(n_2))\tau$, if the individual has
$n_1$ infected neighbours at home and $n_2$ at the workplace. It is
important to note that each hyperedge can be replaced by a clique
(fully connected complete graph) and then the propagation can be
considered on a graph in the conventional way by discounting the
infection pressure for large number of infectious neighbours. In that
case the infection rate will be $f(n_1+n_2)\tau$ instead of
$(f(n_1)+f(n_2))\tau$, hence it cannot be distinguished that somebody
has many infectious neighbours at home and only a few in the
workplace, or a moderate number at both places.

Besides developing the theory of epidemic propagation on hypergraphs
we will run simulations on hypergraphs, derive the exact master
equations and present a mean-field ODE approximation. For the purpose
of simulations we created hypergraphs in several different ways.
%%%
On one hand, the set of vertices was once partitioned into households and once into
workplaces randomly, in this case each node belongs to exactly two
hyperedges. As a second alternative, a Barabási-Albert random graph
was generated and the cliques were determined by a suitable
algorithm. These cliques were considered as the hyperedges of the
newly created hypergraph. This is motivated by the fact that in data mining it is common practice to use different algorithms for identifying cliques as communities from given data of binary relations (i.e. a graph). Finally, we extended the configuration model to generate random hypergraphs with prescribed number and size of hyperedges.

The effect of community structure on the propagation process has been
studied in the literature recently, since non-trivial community
structure occurs not only in epidemiology but also in other systems in
biology, computer science and engineering. Epidemic spreading on
networks with overlapping community structure is considered in
\cite{Chena_etal}, allowing that the rate of infection is different in
the communities. Clique networks describe the phenomenon of
individuals attending different groups in the form of multiple clique
types in \cite{Wang_etal}, where theoretical analysis with a
predefined clique degree distribution is performed. The authors of
\cite{Klamt_etal} conclude that modeling with hypergraphs may result
in a more precise description of biological processes, moreover, they
foresee that "applications of hypergraph theory in computational
biology will increase in the near future".  A voter-like model is
investigated for interacting particle systems on hypergraphs in
\cite{LanchierNeufer}. In their model large blocks of vertices may
flip simultaneously, since vertices in a hyperedge change
simultaneously their opinion to the majority opinion of the hyperedge.

The paper is structured as follows. In Section 2 our model for $SIS$ epidemic propagation on hypergraphs is introduced. Different methods for creating hypergraphs are presented in Section 3. The simulation results for epidemic propagation on different hypergraphs are shown in Section 4. In Section 5 we extend our theory \cite{simon2010JMB} for deriving exact master equations to the case of propagation on hypergraphs. Exact and approximating mean-field equations for the expected value of the number of infected nodes are derived and compared to simulations in Section 6.

%%%%%%%%%%%%%%%%%%%%%%%%%%%%%%%%%%%%%%%%%%%%%%%%%%%%%%%%%%%%%%%%%%%%%%%%%%%%%%%%%%%%%%%%%%
\section{Model formulation}
%%%%%%%%%%%%%%%%%%%%%%%%%%%%%%%%%%%%%%%%%%%%%%%%%%%%%%%%%%%%%%%%%%%%%%%%%%%%%%%%%%%%%%%%%%

Our mathematical model starts from a hypergraph that consists of a set of nodes $V=\{ v_1, v_2, \ldots , v_N\}$ and a set of hyperedges $\mathcal{E}=\{ e_1, e_2, \ldots , e_M\}$ where each hyperedge is a subset of $V$, that is $e_i \subset V$ for all $i=1,2,\ldots ,M$. The pair $(V,\mathcal{E})$ is called a hypergraph. The nodes in our model represent the individuals and the hyperedges are the units of the community structure such as households or workplaces. The methods for creating hypergraphs will be dealt with in the next section, here $(V,\mathcal{E})$ denotes a hypergraph in general. The process considered in this paper is $SIS$ (susceptible-infected-susceptible) epidemic propagation. This means that each node may be in one of the two states susceptible or infected/infectious. These states will be denoted by $S$ and $I$, respectively. A susceptible individual can become infectious after contacting infectious ones and an infectious one can recover and become susceptible again after some time (not depending on the states of its neighbours). Both infection and recovery are governed by a Poisson processes. This means that an infected individual recovers with probability $1-\exp(-\gamma \Delta t)$ in a small time interval $\Delta t$, where $\gamma$ is called the recovery rate. A susceptible individual becomes infected with probability $1-\exp(-r \Delta t)$ in a small time interval $\Delta t$. The rate $r$ is given as
$$
r=\tau \sum_{h} f(k_h) ,
$$
where the summation is for those hyperedges $h\in \mathcal{E}$ that contain the susceptible node, $k_h$ denotes the number of infected nodes in the hyperedge $h$ and $f$ is a given function. We note that in the well-known case of propagation on graphs $f$ is the identity function, i.e. $f(k_h) =k_h$, yielding the infection rate in the form $\tau k$, where $k$ is the total number of infected neighbours. The choice of function $f$ is artificial in this paper. The typical functional form is an inverse tangent like function that is close to the identity around zero and becomes constant for large values of its argument. The behaviour of the propagation process can be tested numerically for a large class of $f$ functions. A careful analysis shows that the qualitative behaviour of the process can be understood by using the simplest piece-wise linear function
\begin{equation}
f(x)= \left\{
\begin{array}{cl}
  x, & \mbox{ if } \quad 0\leq x\leq c \\
  c, & \mbox{ if } \quad x> c.
\end{array} \right. \label{eq:f}
\end{equation}
This function is parametrized by a single parameter $c$, which can be interpreted as a threshold value. If the number of infected nodes in a hyperedge is smaller than $c$, then the infection is proportional to the number of infected neighbours as in the conventional network case, while above this threshold value the number of infected nodes in a hyperedge does not increase the infection pressure. We note that the desirable approach for quantitative analysis would be to determine the functional form and parameters of $f$ by fitting it to real epidemic propagation data. This approach is beyond the scope of this paper and may be the subject of future work.

Thus our model is specified by choosing a hypergraph and the function $f$. Then the full mathematical model is a continuous time Markov chain with a state space of size $2^N$, since to determine the state of the system it has to be given for each node if it is susceptible or infected. The hyperedges of the hypergraph and the function $f$ determine the transition probabilities of the Markov chain as it will be presented in Section \ref{sec:MasterEq}. The full set of master equations of this Markov chain are practically impossible to solve because of the extremely large size of the system, hence individual-based simulations are carried out and the dependence of the system behaviour on the hypergraph structure and on the function $f$ is investigated via simulations, the results of which are shown in Section \ref{sec:Sim}. Before turning to simulation results let us review the methods that were used to create our hypergraphs.

%%%%%%%%%%%%%%%%%%%%%%%%%%%%%%%%%%%%%%%%%%%%%%%%%%%%%%%%%%%%%%%%%%%%%%%%%%%%%%%%%%%%%%%%%%
\section{Hypergraph types used in the simulations}
%%%%%%%%%%%%%%%%%%%%%%%%%%%%%%%%%%%%%%%%%%%%%%%%%%%%%%%%%%%%%%%%%%%%%%%%%%%%%%%%%%%%%%%%%%

As we have mentioned in the introduction  simulations were preformed on hypergraphs generated in three different ways. In this section we introduce these models.

\subsection{Bi-uniform hypergraph model}\label{biuniform}
In the first model, each individual (a vertex of $V$)  belongs to precisely two edges of the hypergraph,
one will correspond to the household
of the individual and the other  to the workplace of the person.
It is assumed that the households are disjoint, and each person works in one workplace.
Thus the set of edges corresponding to the households, $\mathcal{H}$, form a partition of $V$. For sake of simplicity we also assume that each household has precisely $H$ members, thus $(V,\mathcal{H})$ is a $H$-uniform hypergraph consisting of disjoint edges of size $H$. With similar assumptions, the edge set of the workplaces, $\mathcal{W}$, is a partition of $V$ into $W$-element sets. We generate the partitions $\mathcal{H}$ and $\mathcal{W}$ randomly: always pick the next element to the next edge with uniform probability from the remaining unpartitioned elements. Finally we take the union of the two uniform hypergraphs, and obtain a special bi-uniform hypergraph $(V,\mathcal{H}\cup\mathcal{W})$.

\subsection{BA-cliques model}\label{BAcliques}

The second model is based on the preferential attachment model invented by Barabási and Albert \cite{BarabasiAlbert}. They described a randomized algorithm that constructs graphs that are nowadays considered to be one of the best models for numerous natural and human-made structures, including social networks. There are several slightly different ways of implementing this algorithm, however in each case one obtains a graph with similar characteristics.

In \cite{BollobasRiordanST} it was shown that the number of vertices with degree at least $k$ is $k^{-3}$, and in \cite{BollobasRiordan} it was proved that the diameter of such a graph is asymptotically $\log n/\log\log n$. These results match the measures of many natural networks.

In our case, first we generated a random graph using the above method with our own implementation, then it was converted to a hypergraph on the same vertex set. The graph itself represents only the structure of binary relations, however, we need a representation of the structure of some communities, consisting of more members. We assume that each member of a community is in relation with each other member of the community, so members of a community forms a complete subgraph (i.e. there is en edge between any pair of members). Therefore to find these communities we listed all complete subgraphs of our graph using the software called CFinder \cite{Cfinder}, these are the edges of our hypergraph. In this way we obtain many small size hyperedges and few large ones \cite{BianconiMarsili}.

 It is worth mentioning that this algorithm has an exponential worst case running time for general graphs as there may be exponentially many cliques. Fortunately one can show that for graphs obtained by the preferential attachment model this is very unlikely, and the algorithm is most likely efficient in this case.

\subsection{Configuration model}\label{configuration}

\begin{figure}[h!]
\begin{center}
\includegraphics[scale=0.4]{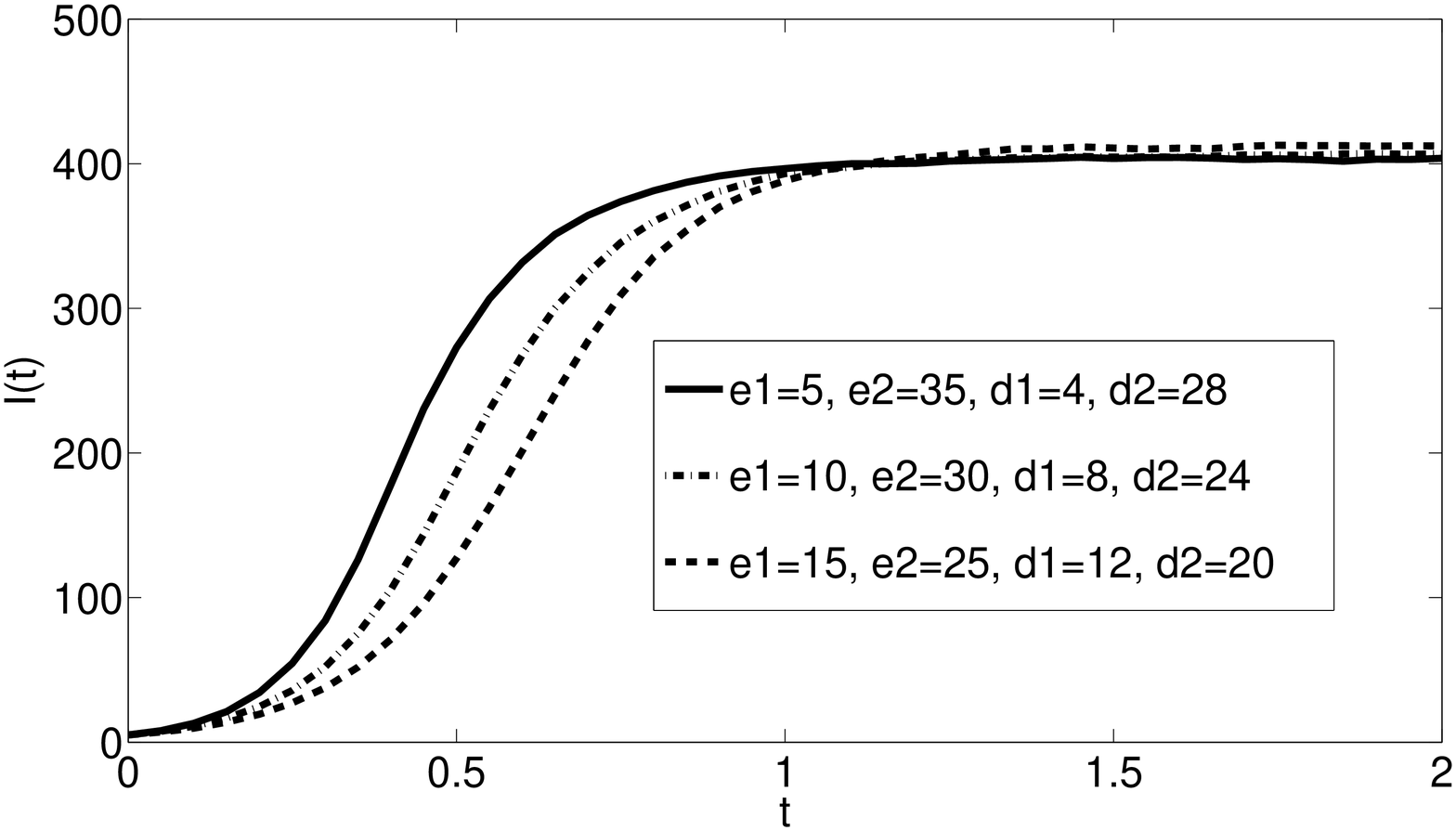}
\caption{Time dependence of the number of infected nodes for bimodal random hypergraphs, in which half of the nodes belong to $d_1$ hyperedges and half of them belong to $d_2$ hyperedges. The total number of hyperedges is $M=400$,  half of them are of size $e_1$, the other half of them are of size $e_2$. The parameter values are $c=10$ in \eqref{eq:f}, $N=500$, $\gamma =1$ and $\tau=0.05$.} \label{fig:random_bimodal}
\end{center}
\end{figure}

In this model our aim is to generate a random, uniform, regular hypergraph on given number of vertices. Therefore, we fix the number of vertices ($N$), the size of each hyperedge ($e$) and the number of edges each vertex is contained in (i.e. the degree of each vertex, $d$). Simple double counting shows that in this way the number of hyperedges is $M=Nd/e$. In a more general setup, one can start with a prescribed size for each hyperedge individually, and a prescribed degree for each vertex. For example, we can set half of the edges to have size $e_1$ and the other half size $e_2$, while half of the vertices have degree $d_1$ and the other half $d_2$ (see Fig. \ref{fig:random_bimodal}.)

All hypergraphs can be associated with a bipartite graph (i.e. a 2-colorable graph) in the following way. Let $V$ be the set of vertices of the hypergraph, this will be one of the color classes in the bipartite graph. The other color class will consist of vertices corresponding to each edge of the hypergraph ($U$). A vertex in $V$ is adjacent to vertex in $U$ if the vertex in contained in edge corresponding to the vertex in $U$. %Let $B(V,\mathcal{E})$ denote this bipartite graph.
Since in our first case every vertex is contained in $d$ hyperedges, each vertex in $V$ has degree $d$, and since each edge contains $e$ vertices, every vertex in $U$ has degree $e$. In the general case, vertices of $V$ have the prescribed degrees, and the degree of the vertices in $U$ will be the prescribed size of the corresponding edge.

Thus it is enough to generate such a random bipartite graph. For this we used the configuration model of Bollob\'as \cite{Bollobasconfig}. For each vertex in the bipartite graph, take the prescribed number of half-edges, then randomly combine two such half edges from the opposite side. In this way parallel edges may occur. It is shown in \cite{Bollobasconfig} that the expected number of parallel edges is relatively small, especially if the prescribed degrees are  small compared to $N$. So we simply chose to delete any parallel edges. As a result a few edges of the corresponding hypergraph contain less than $e$ vertices, and a few vertices are contained in less than $d$ edges. However, it has a negligible effect on the simulation results. The situation is similar in the general case.

%%%%%%%%%%%%%%%%%%%%%%%%%%%%%%%%%%%%%%%%%%%%%%%%%%%%%%%%%%%%%%%%%%%%%%%%%%%%%%%%%%%%%%%%%%
\section{Simulation results} \label{sec:Sim}
%%%%%%%%%%%%%%%%%%%%%%%%%%%%%%%%%%%%%%%%%%%%%%%%%%%%%%%%%%%%%%%%%%%%%%%%%%%%%%%%%%%%%%%%%%

\subsection{Description of the simulation algorithm}

Let $x(t) \in \{ 0,1 \}^N$ denote the state of the system at time $t$. Its $k$-th coordinate, $x_k(t)$, is $0$, if the $k$-th node is susceptible and $1$, if it is infected. The hypergraph is given by its incidence matrix $\mathcal{J}$, the rows of which correspond to the nodes and the coloumns correspond to the hyperedges. That is $\mathcal{J}_{ij}=1$ if node $i$ belongs to the $j$-th hyperedge and it is zero otherwise. Then the product $x(t)\mathcal{J}$ gives the number of infected nodes in the different hyperedges, that is its $j$-th coordinate, $(x(t)\mathcal{J})_j$, is the number of infected nodes in the $j$-th hyperedge.
Since infection and recovery are governed by Poisson processes, a susceptible individual, which has $k_h$ infected neighbours in hyperedge $h$, becomes infected with probability $1-\exp\left(-\tau \sum_{h} f(k_h) \Delta t\right)$ in a small time interval $\Delta t$. Similarly, an infected individual recovers with probability $1-\exp(-\gamma \Delta t)$ in a small time interval $\Delta t$. Thus, assuming that node $i$ is susceptible at time $t$, i.e. $x_i(t)=0$, the rate at which it becomes infected is
$$
\tau \sum_{j=1}^M \mathcal{J}_{ij} f\left((x(t)\mathcal{J})_j \right) .
$$
We apply the widely used individual-base stochastic simulation as for conventional networks. At a given time instant $t$ a vector $r \in [0,1]^N$ is generated with random numbers, then the algorithm runs through all the nodes from $i=1$ to $i=N$. If the $i$-th node is susceptible, i.e. $x_i(t)=0$, then it becomes infected at time $t + \Delta t$, if
\[
r_i < 1-\exp \left(-\tau \sum_{j=1}^M \mathcal{J}_{ij} f\left((x(t)\mathcal{J})_j \right) \Delta t \right).
\]
If the $i$-th node is infected, i.e. $x_i(t)=1$, then it becomes susceptible at time $t + \Delta t$, if
\[
r_i < 1-\exp\left(-\gamma \Delta t\right).
\]
This process is run with sufficiently small time steps $\Delta t$ until the final time $t_{max}$ is reached. We note that running the simulation by using the Gillespie algorithm we obtain completely similar results, when $\Delta t$ is chosen sufficiently small. Then several simulations, started with the same initial condition, are averaged. The simulation results are dealt with in the next subsections.

\subsection{The effect of function $f$}

Here it is studied how the propagation process is affected by the choice of the function $f$. As it was declared in Section 2, in this paper we use the function $f$ given in \eqref{eq:f}. This function is parametrised by a single parameter $c$, which can be interpreted as a threshold value in the following sense. If the number of infected nodes in a hyperedge is smaller than $c$, then the infection is proportional to the number of infected neighbours as in the conventional network case, while above this threshold value the number of infected nodes in a hyperedge do not increase the infection pressure.

\begin{figure}[h!]
\begin{center}
\includegraphics[scale=0.4]{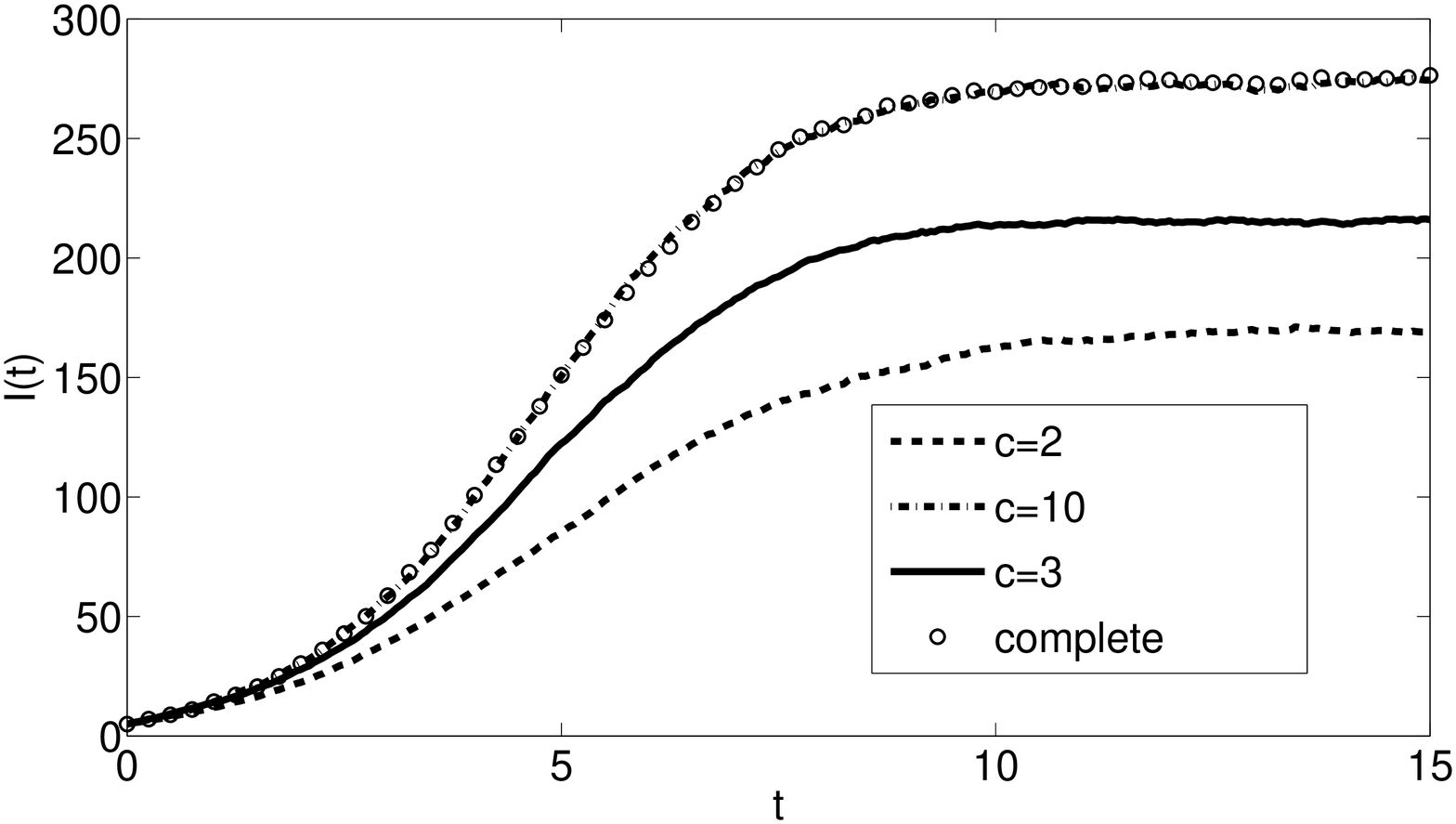}
\caption{Time dependence of the number of infected nodes for hypergraphs with households of size $H=5$ and workplaces of size $W=10$. The curves belonging to different values of the parameter $c$ in \eqref{eq:f} are shown together with the prevalence corresponding to the propagation on a network with hyperedges substituted by complete subgraphs. The parameter values are $N=500$, $\gamma =1$ and $\tau=0.18$.  } \label{fig:f_csalad}
\end{center}
\end{figure}

Consider first the case of bi-uniform hypergraphs constructed from households and workplaces introduced in Subsection~\ref{biuniform}.
A random hypergraph was constructed with households of size $H=5$ and workplaces of size $W=10$, where each node belongs to exactly two hyperedges, a household and a workplace. Then simulations were run with recovery rate $\gamma=1$ and infection rate $\tau=0.18$ for different values of $c$. The time dependence of the number of infected nodes is shown in Figure~\ref{fig:f_csalad} for three different values of $c$. One can see that for greater values of $c$ we face stronger epidemic as it is expected. In order to compare these simulation results to the usual simulation on graphs we created a graph from the hypergraph by exchanging hyperedges to cliques, i.e. to complete subgraphs. In other words, each node in a hyperedge is connected to every other node of that hyperedge with usual edges. It may happen that two nodes are in the same household and in the same workplace, then two edges are created between them in the course of the above process. In this case the two edges are weighted with $2$, in order to make the new network comparable to the original hypergraph. This way a weighted graph is created from the hypergraph. Then simulations were run on this weighted graph as well and the time dependence of the number of infected nodes is compared to that obtained from the simulations on the hypergraph. Figure~\ref{fig:f_csalad} clearly shows that for large values of $c$ the process on the hypergraph is basically the same as the process on the corresponding conventional (weighted) network. This is explained by the simple fact that for $c=10$ the discount effect of the function $f$ cannot come to play. Namely, there are no more than $10$ nodes in the hyperedges, hence $f(k)=k$ for those values of the number of infected neighbours $k$ that can occur in this hypergraph. As it was already mentioned in Section 2, in the course of a quantitative approach this Figure would enable us to fit the value of parameter $c$ by comparing real data to the curves obtained from simulations for different values of $c$.

\begin{figure}[h!]
\begin{center}
\includegraphics[scale=0.4]{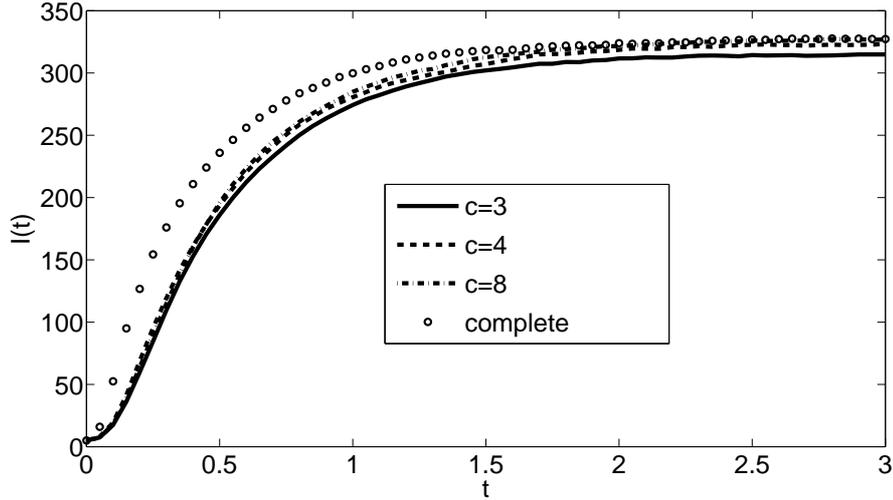}
\caption{Time dependence of the number of infected nodes for hypergraphs created from the cliques of a Barab\'asi-Albert graph. The curves belonging to different values of the parameter $c$ in \eqref{eq:f} are shown together with the prevalence corresponding to the propagation on a network with hyperedges substituted by complete subgraphs. The parameter values are $N=500$, $\gamma=1$ and $\tau=0.02$ .} \label{fig:barabasi}
\end{center}
\end{figure}

Consider now the case of hypergraphs created from the cliques of Barabási-Albert networks, see Subsection~\ref{BAcliques}. We constructed a network with $N=500$ nodes by using the preferential attachment model with $m=4$. Then the cliques were determined and substituted by hyperedges. Simulations were run with recovery rate $\gamma=1$ and infection rate $\tau=0.02$ for different values of $c$. The time dependence of the number of infected nodes is shown in Figure~\ref{fig:barabasi} for three values of $c$. One can again see that for greater values of $c$ we face stronger epidemic as it is expected. The smallest  value, $c=3$, is below the average hyperedge size, hence the effect of the function $f$ can be clearly seen. Namely, the infection is smaller than on the weighted graph, which is created from the hypergraph by changing hyperedges to cliques. On the other hand, for the largest value, $c=8$, which is greater than the average size of the hyperedges, the hypergraph structure has negligible effect in the steady state compared to the conventional propagation on the weighted network. More importantly, we can see that the hyperedge model has strong effect in the early stage of the epidemic, when large cliques are infected probably, for which the hyperedge size is larger than these values of $c$.

\begin{figure}[h!]
\begin{center}
\includegraphics[scale=0.4]{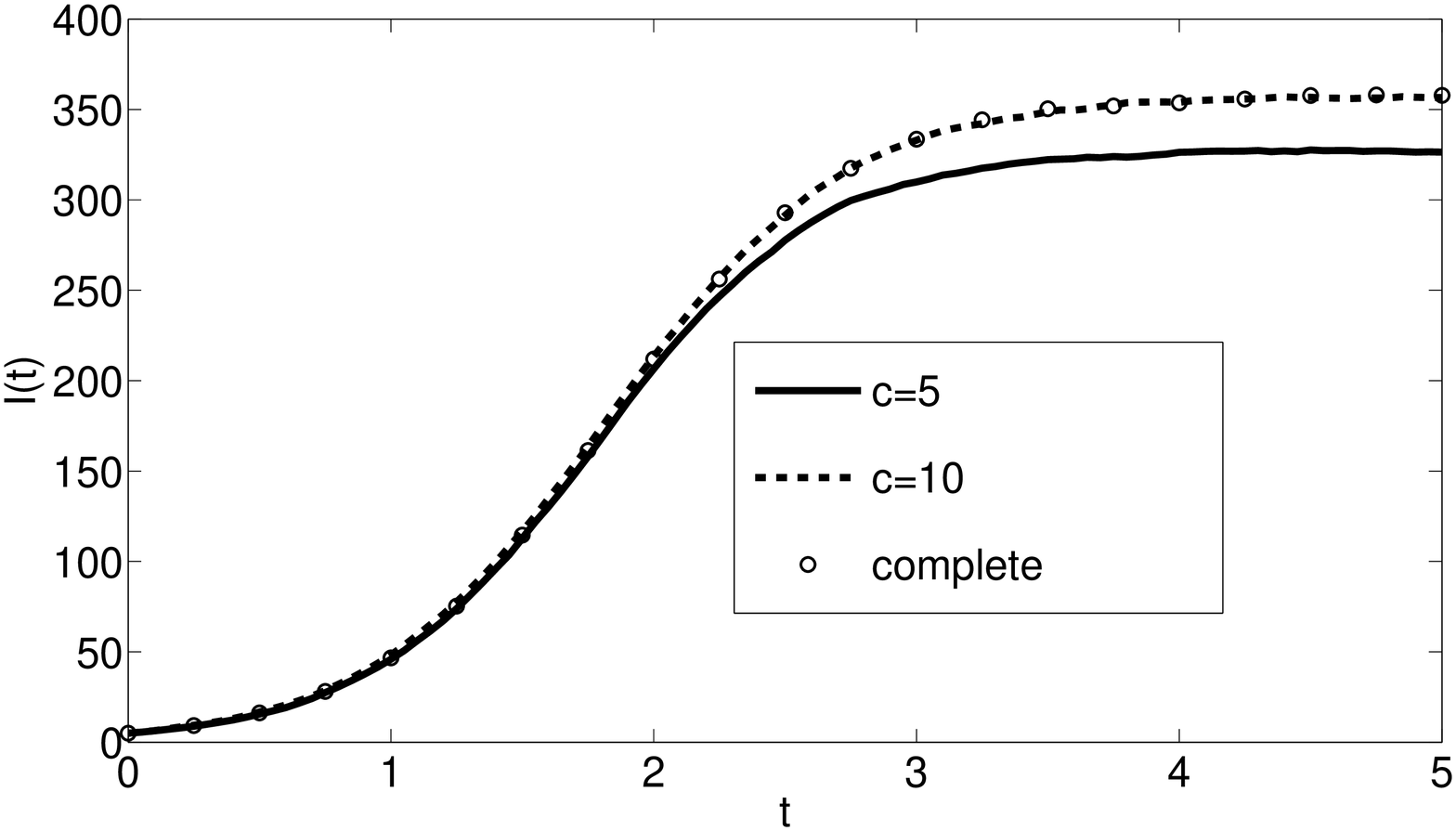}
\caption{Time dependence of the number of infected nodes for a regular random hypergraph, in which every node belongs to $d=8$ hyperedges, each of which is of size $e=10$. The curves belonging to different values of the parameter $c$ in \eqref{eq:f} are shown together with the prevalence corresponding to the propagation on a network with hyperedges substituted by complete subgraphs. The parameter values are $N=500$, $M=400$, $\gamma =1$ and $\tau=0.05$.} \label{fig:f_random}
\end{center}
\end{figure}

Consider now the case of random hypergraphs created by the configuration model, see Subsection~\ref{configuration}. We constructed a regular hypergraph, in which all hyperedges contain $e=10$ nodes and the degree of each node is $d=8$, i.e. each node belongs to $8$ hyperedges. Then simulations were run with recovery rate $\gamma=1$ and infection rate $\tau=0.05$ for different values of $c$. The time dependence of the number of infected nodes is shown in Figure~\ref{fig:f_random} for two values of $c$. One can again see that for greater values of $c$ we face stronger epidemic as it is expected. The smaller value, $c=5$, is below the hyperedge size, hence the effect of the function $f$ can be observed. Namely, the infection is smaller than on the weighted graph, which is created from the hypergraph by changing hyperedges to cliques. On the other hand, for the greater value $c=10$, which is the same as the size of the hyperedges, the hypergraph structure has no effect compared to the conventional propagation on  network, since $f(k)=k$ for those values of the number of infected neighbours $k$ that can occur in this hypergraph.

\subsection{The effect of the structure of the hypergraph}

\begin{figure}[h!]
\begin{center}
\includegraphics[scale=0.4]{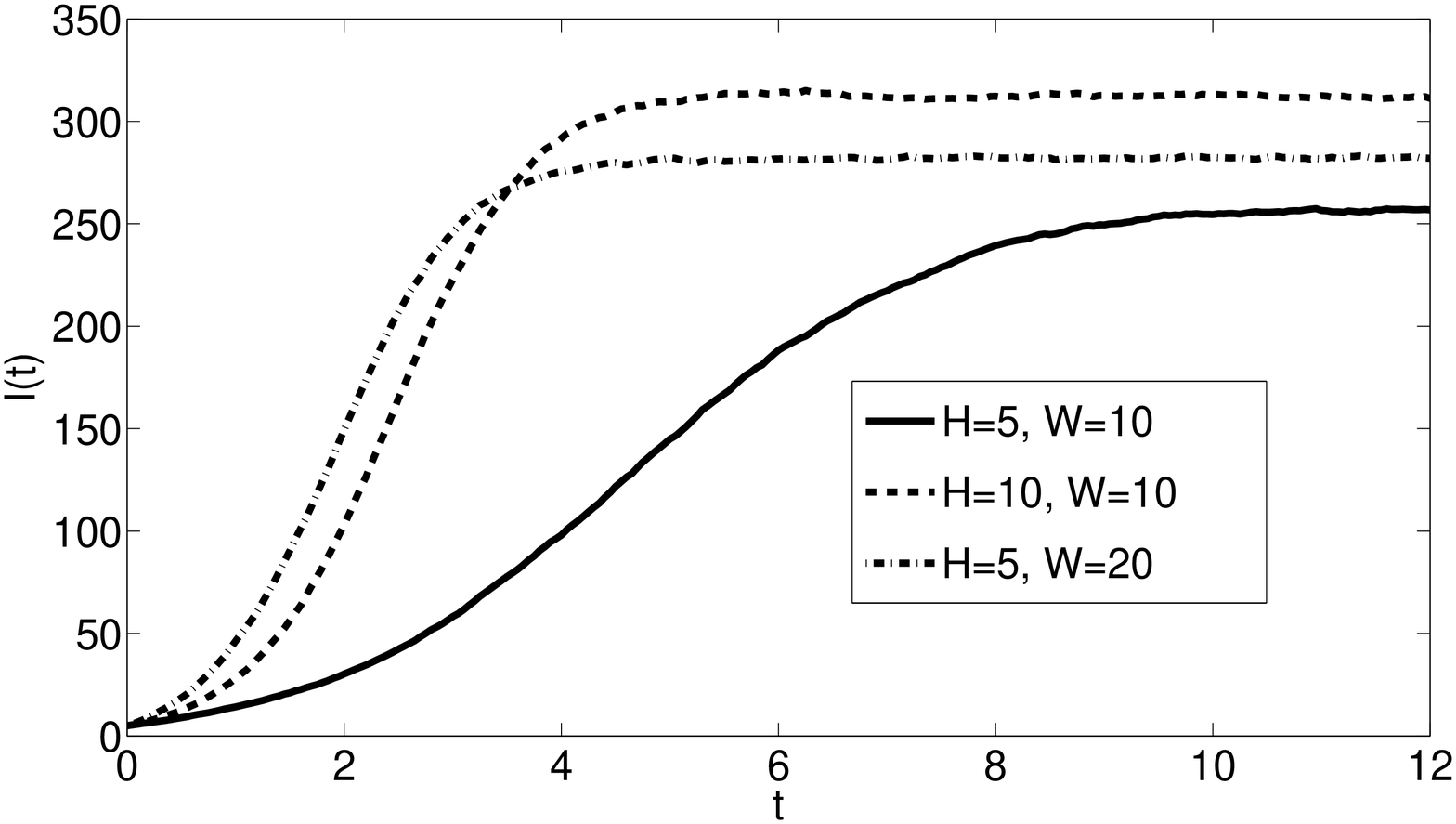}
\caption{Time dependence of the number of infected nodes for hypergraphs with households of size $H$ and workplaces of size $W$. The parameter values are $c=5$ in \eqref{eq:f}, $N=500$, $\gamma =1$ and $\tau=0.18$.} \label{fig:hyper_HW}
\end{center}
\end{figure}

Now we fix the function $f$, that is fix a value of $c$ and investigate how the parameters of the hypergraph affect the spreading process. Consider first the case of bi-uniform hypergraphs constructed from households and workplaces in Subsection~\ref{biuniform}. A random hypergraph was constructed with households of size $H$ and workplaces of size $W$, where each node belongs to exactly two hyperedges, a household and a workplace. Then simulations were run with recovery rate $\gamma=1$ and infection rate $\tau=0.18$ for different values of $H$ and $W$. Figure~\ref{fig:hyper_HW} shows that the number of infected nodes during the spreading process depends in a non-trivial way on the values of $H$ and $W$. One can observe that increasing the size of the hyperedges, for example, comparing the case $H=5$, $W=10$ to the case $H=10$, $W=10$, the infection becomes stronger. On the other hand, comparing the case $H=5$, $W=20$ to the case $H=10$, $W=10$ we can say that on a more heterogeneous (in the sense of hyperedge sizes) hypergraph the propagation starts faster, but it can settle at a smaller steady state value than in the case of a more homogeneous hypergraph.

Let us turn to the study of the effect of degree heterogeneity. A hypergraph is constructed with $M$ hyperedges, half of them is of size $e_1$ and and the other half is of size $e_2$. Half of the $N$ nodes has degree $d_1$, i.e. belong to $d_1$ hyperedges, and the half of them is of degree $d_2$. We note that the conservation relation $\frac{M}{2} (e_1 +e_2)= \frac{N}{2} (d_1 +d_2)$ must hold. In Figure~\ref{fig:random_bimodal} the time dependence of the number of infected nodes is shown for different degree distributions. We can see that for a more homogeneous degree distribution the spreading starts slower but it ends at a higher value. It is also important to note that for the most heterogeneous case, when the size $e_1$ of the small hyperedges is less than $c$, the discount effect of the function $f$ applies only for the large hyperedges, while in the most homogeneous case with $e_1=15$, $e_2=25$ the function $f$ affects all hyperedges.

\begin{figure}[h!]
\begin{center}
\includegraphics[scale=0.4]{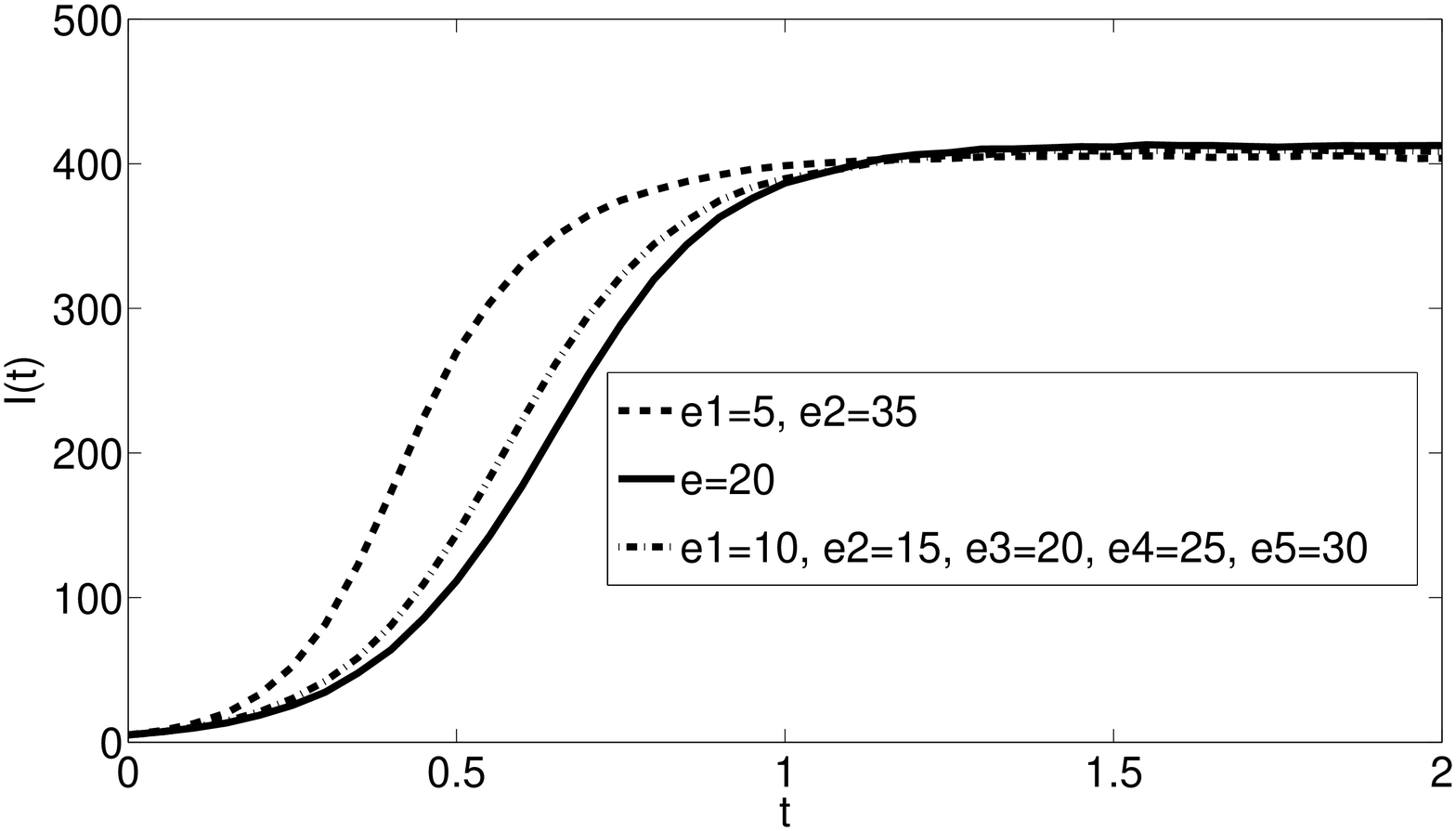}
\caption{Time dependence of the number of infected nodes for a regular hypergraph (continuous line), for a bimodal hypergraph (dashed line) and for a hypergraph with 5 different hyperedge sizes (dashed-dotted line), with the same number of hyperedges in each category. The parameter values are $c=10$ in \eqref{eq:f}, $N=500$, $M=400$,   $\gamma =1$ and $\tau=0.05$.}\label{fig:random_multi_e}
\end{center}
\end{figure}

The effect of heterogeneity is also shown in Figure~\ref{fig:random_multi_e}, where the spread on a regular and a bimodal hypergraph is compared to that on a hypergraph with five different hyperedge sizes. The hypergraph have $N=500$ nodes and consist of $M=400$ hyperedges. The number of hyperedges is the same in each size category, i.e. for the bimodal hypergraph there are 200 hyperedges with both sizes $e_1$ and $e_2$ and for the hypergraph with five hyperedge sizes there are 80 hyperedges with sizes $e_i$ ($i=1,2,\ldots , 5$). One can again observe that the fastest spread is on the most heterogeneous hypergraph in the initial phase of the process, while this network leads to the least prevalence in the final stage of the process.

%%%%%%%%%%%%%%%%%%%%%%%%%%%%%%%%%%%%%%%%%%%%%%%%%%%%%%%%%%%%%%%%%%%%%%%%%%%%%%%%%%%%%%%%%%
\section{Master equation} \label{sec:MasterEq}
%%%%%%%%%%%%%%%%%%%%%%%%%%%%%%%%%%%%%%%%%%%%%%%%%%%%%%%%%%%%%%%%%%%%%%%%%%%%%%%%%%%%%%%%%%

Here the exact master equations of SIS epidemic propagation on an arbitrary hypergraph with $N$ nodes are derived. The master equations form a linear system of ordinary differential equations, the coefficients of which are the transition rates from one state to another.

The state space of the process is $\{ S,I \}^N$ containing $2^N$ elements, since each node can be in one of two states $S$ or $I$. The state space is divided into $N+1$ classes according to the number of infected nodes. Let $S^0$ denote the state where every node is susceptible, i.e. $S^0=(S,S,\ldots,S)$.
Let $S^k$ be the subset with states having $k$ infected nodes, containing $c_k= \binom{n}{k}$ states. Finally, let $S^N$ be the state in which every node is infected, i. e. $S^N=(I,I,\ldots,I)$.

The elements of $S^k$ are denoted by $S_1^k, S_2^k, \ldots, S_{c_k}^{k}$. Let $S_j^k(l)$ be the type of the $l$-th node in the state $S_j^k$, so that the value of $S_j^k(l)$ is either $S$ or $I$. The state of the system can change in two ways:

\begin{enumerate}
\item Infection: a susceptible node becomes infected, this is an $S_j^k \to S_i^{k+1}$ transition, where $i$ and $j$ are such that there exists $l$, for which $S_j^k(l)=S$, $S_i^{k+1}(l)=I$ and $S_j^k(m)=S_i^{k+1}(m)$ for every $m \neq l$. Furthermore, node $l$ has an infected neighbour, which can be formally expressed as follows: there exists an $r \neq l$, such that $S_j^k(r)=I$ and $l$ and $r$ are in the same hyperedge.

\item Recovery: an infected node becomes susceptible, this is an $S_j^k \to S_i^{k-1}$ transition, where $i$ and $j$ are such that there exists $l$, for which $S_j^k(l)=I$, $S_i^{k-1}(l)=S$ and $S_j^k(m)=S_i^{k-1}(m)$ for every $m \neq l$.
\end{enumerate}

Let $X_j^k(t)$ denote the probability that the system is in state $S_j^k$ at time $t$.
Let
\begin{equation*}
X^k(t)=\left(X_1^k(t),X_2^k(t), \ldots, X_{c_k}^k(t)\right)
\end{equation*}
denote the probability of states containing $k$ infected nodes, where $k=0,1,\ldots,N$.
The above transitions define a system of linear differential equations with constant coefficients for $X_j^k(t)$, known as the Kolmogorov equations or master-equations.
The number of the infectious nodes changes by one at most in each time step, thus the master-equation can be written in the following block tridiagonal form:
\begin{equation}\label{(5.1)}
\dot{X^k}=A^k X^{k-1} + B^k X^k + C^k X^{k+1}, \qquad k=0,1,\ldots, N,
\end{equation}
where $A^0$ and $C^N$ are zero matrices.
In matrix form:
\begin{equation*}
\dot{X}=P X,
\end{equation*}
where
\[
P =\left(
\begin{array}{cccccc}
  B^0 & C^0 & 0 & 0 & 0 & 0 \\
  A^1 & B^1 & C^1 & 0 & 0 & 0\\
  0 & A^2 & B^2 & C^2 & 0 & 0 \\
  0 & 0 & A^3 & B^3 & C^3 & 0 \\
  \vdots & \vdots & \cdots & \cdots & \cdots & \vdots \\
  0 & 0 & \cdots & \cdots & A^N & B^N
 \end{array} \right).
\]
The matrices $A^k$ describe infection and $C^k$ describe recovery. The structure of the network is reflected by the matrices $A^k$.

Let $A_{i,j}^k$ denote the $(i,j)$-th element of $A^k$, which shows the transition rate from state $S_j^{k-1}$ to state $S_i^k$.
The class $S^{k-1}$ has $c_{k-1}$ elements and the class $S^{k}$ consists of $c_k$ terms, therefore the matrix $A^k$ has $c_k$ rows and $c_{k-1}$ columns.
The entry $A_{i,j}^k$ is non-zero if and only if $S_j^{k-1}$ and $S_i^k$ differ only at one position. Let the $l$-th node be this one, i.e. $S_j^{k-1}(l)=S$, $S_i^k(l)=I$ and $S_j^{k-1}(m)=S_i^k(m)$ for every $m \neq l$. Furthermore, there exists a number $r \neq l$, such that $S_j^{k-1}(r)=I$ and $l$ is in the same hyperedge as $r$.
Then the transition rate is given as
\begin{equation}\label{(5.2)}
A_{i,j}^k = \tau \sum_{h:\, l \in h} f\left(N_h(S_j^{k-1})\right),
\end{equation}
where $N_h(S_j^{k-1})$ denotes the number of infected nodes in hyperedge $h$ in the state $S_j^{k-1}$ and $f$ is the function given in \eqref{eq:f}. Summing these equations for $i$ one obtains for every $j \in \{1,2,\ldots,c_{k-1}\}$ that
\begin{equation}\label{(5.3)}
\sum_{i=1}^{c_k} A_{i,j}^k = \tau N_{SI}^f(S_j^{k-1}),
\end{equation}
where
\begin{equation}\label{NSIf}
N_{SI}^f(S_j^k) = \sum_{l:\, S_j^k(l)=S} \,\, \sum_{h:\, l \in h} f\left(N_h(S_j^k)\right).
\end{equation}
That is $N_{SI}^f(S_j^k)$ denotes the sum of the values $\sum_{h: \, l \in h} f(N_h(S_j^{k+1}))$ for susceptible nodes in state $S_j^k$.

Let $C_{i,j}^k$ denote the element in the $i$-th row and $j$-th column of the matrix $C^k$.
This gives the transition rate from state $S_j^{k+1}$ to state $S_i^k$.
The class $S^{k+1}$ has $c_{k+1}$ elements and the class $S^{k}$ has $c_k$ terms, therefore $C^k$ has $c_k$ rows and $c_{k+1}$ columns.
The entry $C_{i,j}^k$ is non-zero if and only if $S_j^{k+1}$ and $S_i^k$ differ only at one node. Let this node be the $l$-th one, that is $S_j^{k+1}(l)=I$, $S_i^k(l)=S$ and $S_j^{k+1}(m)=S_i^k(m)$ for every $m \neq l$.
In this case $C_{i,j}^k=\gamma$. The number of infected nodes in state $S_j^{k+1}$ is $k+1$, hence there are $k+1$ elements in the $j$-th column of the matrix $C^k$, which are equal to $\gamma$, all other entries are zero. Thus for every $j \in \{1,2,\ldots,c_{k+1} \}$ we have
\begin{equation}\label{(5.4)}
\sum_{i=1}^{c_k} C_{i,j}^k=\gamma (k+1).
\end{equation}

The $B^k$ matrices are diagonal with $c_k$ rows and coloumns. The elements of $B^k$ denote the rate of the $S_i^k \to S_j^k$ type transition, which are non-zero if and only if $i = j$. Since the column-wise sum of the elements of $P$ are zero, the matrix $B^k$ is determined as follows:
\begin{equation}\label{(5.5)}
B_{i,i}^k = -\sum_{j=1}^{c_{k+1}} A_{j,i}^{k+1} - \sum_{j=1}^{c_{k-1}} C_{j,i}^{k-1}.
\end{equation}

\begin{figure}[ht!]
\centering
\begin{tikzpicture}[scale=1]
\tikzstyle{every node}=[circle]
\begin{scope}[fill opacity=0.5]
\filldraw[fill=gray!400] (2,0) ellipse (3.5cm and 1cm);
\filldraw[fill=gray] (4,2) ellipse (1cm and 3.5cm);
\filldraw[rounded corners=1cm, fill=gray!30] (-1,-1)--(1,-1)--(1,3)--(5,3)--(5,5)--(-1,5)--cycle;
\end{scope}
\node at (0,0) [label={4},draw,fill] {};
\node at (0,4) [label={1},draw,fill] {};
\node at (4,4) [label={2},draw,fill] {};
\node at (4,0) [label={3},draw,fill] {};
\end{tikzpicture}
\caption{A simple hypergraph with 4 nodes and 3 hyperedges.} \label{fig:hypgraph}
\end{figure}
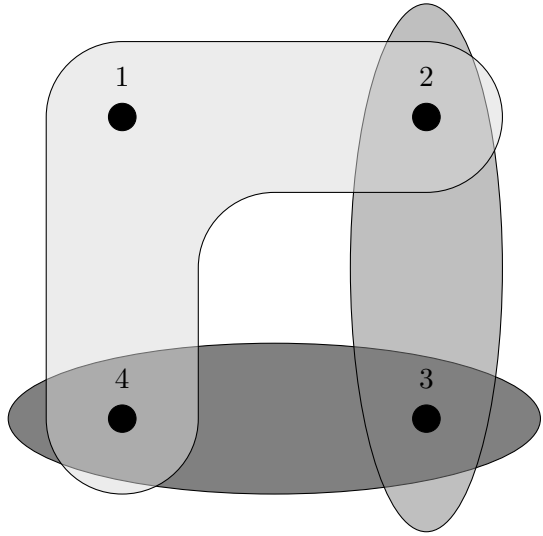

As an example, the master equations are determined below for the hypergraph with $N=4$ nodes shown in Figure \ref{fig:hypgraph}. The incidence matrix of the hypergraph is
\[
\mathcal{J} =\left(
\begin{array}{ccc}
  1 & 0 & 0 \\
  1 & 1 & 0\\
  0 & 1 & 1\\
  1 & 0 & 1
 \end{array} \right).
\]
In this case the state space has $2^4$ elements and we have four classes according to the number of infected nodes:
\begin{align*}
X^0 &= X_{SSSS}, \\
X^1 &= \left(X_{SSSI}, X_{SSIS}, X_{SISS}, X_{ISSS}\right),\\
X^2 &= \left(X_{SSII}, X_{SISI}, X_{SIIS}, X_{ISSI}, X_{ISIS}, X_{IISS}\right),\\
X^3 &= \left(X_{SIII}, X_{ISII}, X_{IISI}, X_{IIIS}\right),\\
X^4 &= X_{IIII} .
\end{align*}
The vector of probabilities is $X=\left(X^0, X^1, X^2, X^3, X^4\right)$ and the matrix $P$ of the linear system of ODEs yielding the master equations is
\[
P =\left(
\begin{array}{ccccc}
  B^0 & C^0 & 0 & 0 & 0 \\
  A^1 & B^1 & C^1 & 0 & 0 \\
  0 & A^2 & B^2 & C^2 & 0 \\
  0 & 0 & A^3 & B^3 & C^3 \\
  0 & 0 & 0 & A^4 & B^4
 \end{array} \right).
\]
The submatrices are be given as follows
\[
A^1 =\left(
\begin{array}{c}
  0 \\
  0 \\
  0 \\
  0
 \end{array} \right),
\qquad
A^2 = \tau \left(
\begin{array}{cccc}
  f(1) & f(1) & 0 & 0 \\
  f(1) & 0 & f(1) & 0 \\
  0 & f(1) & f(1) & 0 \\
  f(1) & 0 & 0 & f(1) \\
  0 & 0 & 0 & 0 \\
  0 & 0 & f(1) & f(1)
 \end{array} \right),
\]
\[
A^3 = \tau \left(
\begin{array}{cccccc}
  2 f(1) & 2 f(1) & 2 f(1) & 0 & 0 & 0 \\
  f(1) & 0 & 0 & f(1) & 2 f(1) & 0 \\
  0 & f(2) & 0 & f(2) & 0 & f(2) \\
  0 & 0 & f(1) & 0 & 2 f(1) & f(1) \\
 \end{array} \right),
\]
\[
A^4 =\tau \left( f(2), f(1)+f(2), 2 f(1), f(1) + f(2) \right),
%\left(
%\begin{array}{cccc}
%  f(2) & (f(1)+f(2)) & 2 f(1) & (f(1)+f(2))
% \end{array} \right).
\]
\[
B^0 = ( 0 ),
%\left(
%\begin{array}{c}
%  0
% \end{array} \right),
\qquad
C^0 = \left( \gamma, \gamma, \gamma, \gamma \right),
%\left(
%\begin{array}{cccc}
%  \gamma & \gamma & \gamma & \gamma
% \end{array} \right),
\]
\[ C^1 = \left(
\begin{array}{cccccc}
\gamma & \gamma & 0 & \gamma & 0 & 0\\
\gamma & 0 & \gamma & 0 & \gamma & 0 \\
0 & \gamma & \gamma & 0 & 0 & \gamma \\
0 & 0 & 0 & \gamma & \gamma & \gamma \\
\end{array} \right),
\]
\[
C^2 = \left(
\begin{array}{cccc}
\gamma & \gamma & 0 & 0 \\
\gamma & 0 & \gamma & 0 \\
\gamma & 0 & 0 & \gamma \\
0 & \gamma & \gamma & 0 \\
0 & \gamma & 0 & \gamma \\
0 & 0 & \gamma & \gamma
\end{array} \right), \qquad
C^3 = \left(
\begin{array}{c}
\gamma \\
\gamma \\
\gamma \\
\gamma
\end{array} \right).
\]

For example, the elements of the third row of the matrix $A^3$ are
\[
\left(0, \tau f(2), 0, \tau f(2), 0, \tau f(2) \right),
\]
which describe the following rates of transitions: $SSII \to IISI$, $SISI \to IISI$, $SIIS \to IISI$, $ISSI \to IISI$, $ISIS \to IISI$, $IISS \to IISI$.
The first, third and fifth transitions can not be realized, as the states differ in more than one node. During the second transition, the first node becomes infected, hence we need to compute the term $N_{SI}^f(SISI)$. This is equal to $f(2)$, because the first node is in a hyperedge, where it has two infected  neighbours. The other transitions can be determined in a similar way.

Using this method one can determine the master equations for an arbitrary hypergraph theoretically. However, we have to note that the master equations are useful mainly from the theoretical point of view, as the above construction can only be carried out practically for hypergraphs of moderate size, because of the huge number of equations. This motivates the derivation of approximating systems, called mean-field equations that will be dealt with in the next section.

%%%%%%%%%%%%%%%%%%%%%%%%%%%%%%%%%%%%%%%%%%%%%%%%%%%%%%%%%%%%%%%%%%%%%%%%%%%%%%%%%%%%%%%%%%
\section{Mean-field theory}
%%%%%%%%%%%%%%%%%%%%%%%%%%%%%%%%%%%%%%%%%%%%%%%%%%%%%%%%%%%%%%%%%%%%%%%%%%%%%%%%%%%%%%%%%%

\subsection{Exact differential equations for the expected number of infected and susceptible nodes}

The main idea of mean-field theory is to consider some expected quantities instead of the probabilities of each individual state. The most important quantities are the expected number of infected and susceptible nodes that can be given as
\begin{equation}\label{(5.11)}
[I](t) = \sum_{k=0}^N k \sum_{j=1}^{c_k} X_j^k(t), \qquad [S](t) = \sum_{k=0}^N (N-k) \sum_{j=1}^{c_k} X_j^k(t).
\end{equation}
During the derivation of differential equations for these expected values below, we will need the quantity
\begin{equation}\label{(5.12)}
[SI]=\sum_{k=0}^N \sum_{j=1}^{c_k} N_{SI}^f(S_j^k) X_j^k(t),
\end{equation}
that can be considered as the generalization of the average number of $SI$ edges defined in conventional networks. We remind that $N_{SI}^f$ is defined in \eqref{NSIf}. Now we are in the position to derive the exact differential equations for the expected values $[I](t)$ and $[S](t)$ starting from the master equations \eqref{(5.1)}.

\begin{theo} \label{theo:MF}
The expected values $[I](t)$ and $[S](t)$ satisfy the following differential equations for an arbitrary hypergraph.
\begin{align}
\dot{[S]} &= \gamma [I] - \tau [SI],\label{(5.13)} \\
\dot{[I]} &= \tau [SI] - \gamma [I]. \label{(5.14)}
\end{align}
\end{theo}

%\begin{proof}

\noindent \textit{Proof}

Introducing the notation $S_k = (1, \, 1, \, \ldots , 1)$ we have
\begin{equation*}
\sum_{j=1}^{c_k} X_j^k = S_k X^k
\end{equation*}
hence \eqref{(5.11)} takes the form
\begin{equation}\label{(5.15)}
[I](t)=\sum_{k=0}^N k S_k X^k, \qquad [S](t)=\sum_{k=0}^N (N-k) S_k X^k.
\end{equation}
Equation \eqref{(5.5)} can be written as
\begin{equation*}
B_{i,i}^k = -\left(S_{k+1} A^{k+1}\right)_i - \left(S_{k-1} C^{k-1}\right)_i,
\end{equation*}
and using that $B^k$ is a diagonal matrix we get
\begin{equation*}
B_{i,i}^k = \left(S_k B^k \right)_i.
\end{equation*}
Thus for every $i=1,\ldots,c_k$
\begin{equation*}
\left( S_k B^k \right)_i = - \left( S_{k+1} A^{k+1} \right)_i - \left( S_{k-1} C^{k-1} \right)_i
\end{equation*}
holds that can be written as
\begin{equation}\label{(5.16)}
S_{k+1} A^{k+1} + S_k B^k + S_{k-1} C^{k-1} = 0,
\end{equation}
holding for every $k=0,1, \ldots, N$, where $A^{N+1}$ and $C^{-1}$ are zero matrices.

Differentiating the function $[I](t)$ and using the derivatives of $X^k$ given by \eqref{(5.1)} leads to
\begin{gather*}
\dot{[I]}=\sum_{k=0}^N k S_k \dot{X^k} = \sum_{k=0}^N k S_k \left( A^k X^{k-1} + B^k X^k + C^k X^{k+1} \right) = \\
= \sum_{k=1}^N k S_k A^k X^{k-1} + \sum_{k=0}^N k S_k B^k X^k + \sum_{k=0}^{N-1} k S_k C^k X^{k+1} = \\
= \sum_{k=0}^{N-1} (k+1) S_{k+1} A^{k+1} X^k + \sum_{k=0}^N k S_k B^k X^k + \sum_{k=1}^N (k-1) S_{k-1} C^{k-1} X^k = \\
= \sum_{k=0}^N \left( (k+1) S_{k+1} A^{k+1} + k S_k B^k + (k-1) S_{k-1} C^{k-1} \right) X^k.
\end{gather*}
Thus from equation \eqref{(5.16)} we obtain the following differential equation:
\begin{equation*}
\dot{[I]}=\sum_{k=0}^N \left( S_{k+1} A^{k+1} - S_{k-1} C^{k-1} \right) X^k.
\end{equation*}
Now, the desired equation, \eqref{(5.14)} can be derived by using the proposition below. The proof for $[S](t)$ is similar.
%\end{proof}

\begin{prop}
The matrices $A^k$ and $C^k$ satisfy the following identities.
\begin{enumerate}
\item $S_{k-1} C^{k-1} = \gamma k S_k,$
\item $\sum\limits_{k=0}^N S_{k-1} C^{k-1} X^k = \gamma [I],$
\item $\sum\limits_{k=0}^{N} S_{k+1} A^{k+1} X^k = \tau [SI].$
\end{enumerate}
\end{prop}

%\begin{proof}
\noindent \textit{Proof}

From equation \eqref{(5.4)} one obtains
\begin{equation*}
\left( S_{k-1} C^{k-1} \right)_j = \sum_{i=1}^{c_{k-1}} C_{i,j}^{k-1} = \gamma k,
\end{equation*}
for every $j \in \{ 1,2,\ldots,c_k \}$, therefore $S_{k-1} C^{k-1} = \gamma k S_k$, proving the first part of the statement.

The second part follows from the first one by using equation \eqref{(5.15)}.

To prove the last statement write equation \eqref{(5.3)} as
\begin{equation*}
\left( S_{k+1} A^{k+1} \right)_j = \sum_{i=1}^{c_{k+1}} A_{i,j}^{k+1} = \tau N_{SI}^f(S_j^k),
\end{equation*}
for an arbitrary $j \in \{ 1,2,\ldots, c_k \}$. Therefore
\begin{equation*}
\sum_{k=0}^N S_{k+1} A^{k+1} X^k = \sum_{k=0}^N \sum_{j=1}^{c_k} \left( S_{k+1} A^{k+1} \right)_j X_j^k = \tau \sum_{k=0}^N \sum_{j=1}^{c_k} N_{SI}^f(S_j^k) X_j^k(t) = \tau [SI]
\end{equation*}
that we wanted to prove.
%\end{proof}

Theorem \ref{theo:MF} provides exact differential equations for the expected number of susceptible and infected nodes, however, these differential equations are not self-contained, since the function $[SI]$ is not known. The next step of the mean-field theory is to derive an approximation of $[SI]$ in terms of $[S]$ and $[I]$ in order to make the system closed. This is called a moment closure approximation, which is the subject of the next subsection.

\subsection{Closed mean-field equations and their comparison to simulation} \label{sec:MF}

Our aim now is to derive an approximation to $[SI]$ given in \eqref{(5.12)}, which is based on the definition in \eqref{NSIf}.

\begin{figure}[h!]
\begin{center}
\includegraphics[scale=0.5]{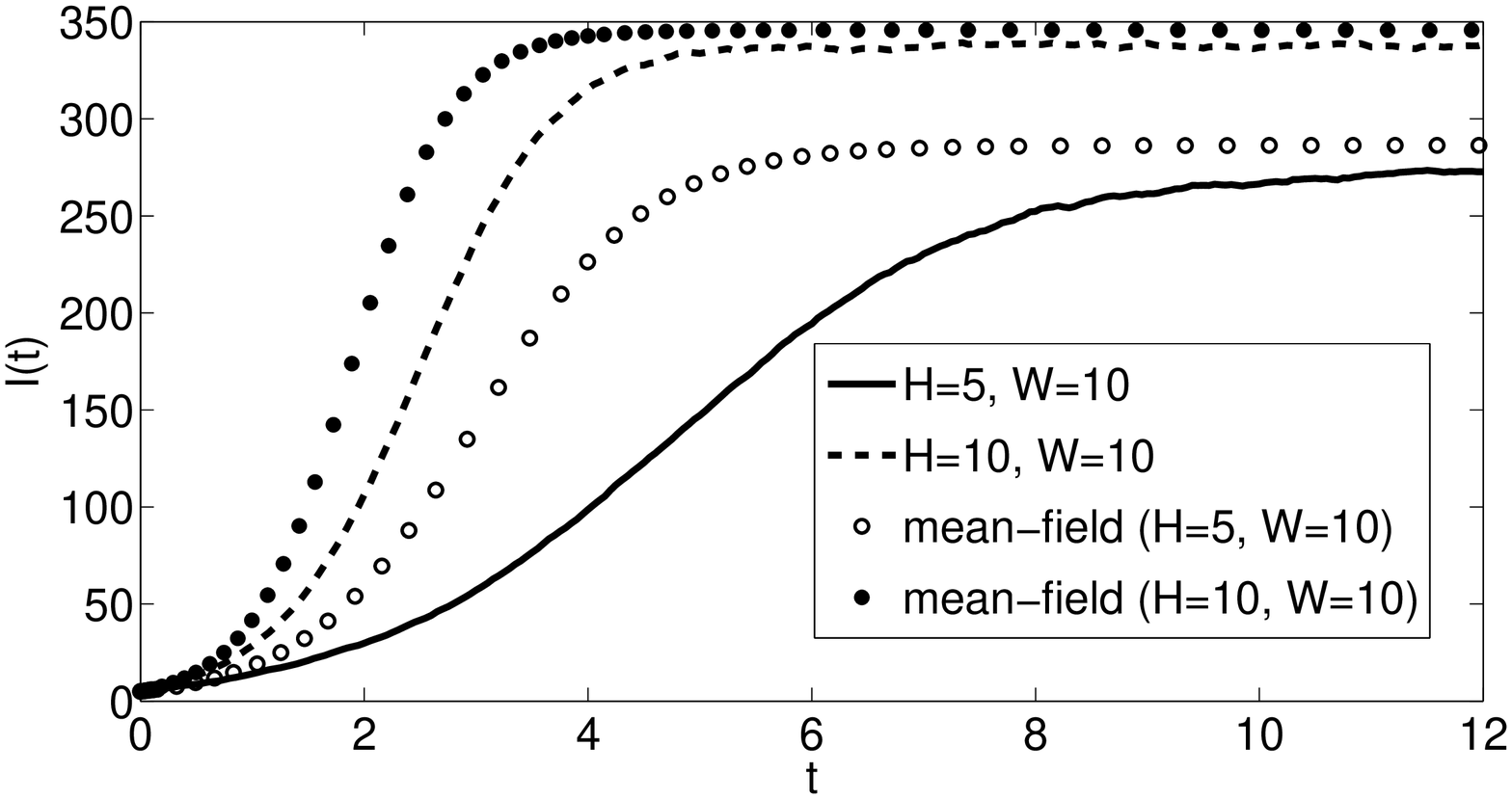}
\caption{Time dependence of the number of infected nodes for hypergraphs with households of size $H$ and workplaces of size $W$. The solution of the mean-field equation \eqref{csaladode45} is also shown. The parameter values are $c=7$ in \eqref{eq:f}, $N=500$, $\gamma =1$ and $\tau=0.18$. } \label{fig:csaladok}
\end{center}
\end{figure}

Consider first hypergraphs describing the network structure with households of size $H$ and workplaces of size $W$. In these hypergraphs, a node belongs to exactly two hyperedges, a household and a workplace. Denoting the total number of infected nodes by $I$, the average number of infected neighbours of a node can be approximated by $\frac{H-1}{N} I$ in a household and by $\frac{W-1}{N}I$ in a workplace. Thus the second summation in \eqref{NSIf} consists of two terms and is approximated as
$$
f\left(\frac{H-1}{N} I \right)+f\left(\frac{W-1}{N} I \right) .
$$
Since this is independent of $l$, the double sum in \eqref{NSIf} reduces to
$$
N_{SI}^f(S_j^k) \approx (N-k) \left[ f\left(\frac{H-1}{N} I \right)+f\left(\frac{W-1}{N} I \right) \right],
$$
because the number of susceptible nodes in state $S_j^k$ is $N-k$. Now \eqref{(5.12)} leads to the approximation
$$
[SI]\approx \left[ f\left(\frac{H-1}{N} I \right)+f\left(\frac{W-1}{N} I \right) \right] \sum_{k=0}^N \sum_{j=1}^{c_k} (N-k) X_j^k(t) =
$$
$$
\left[ f\left(\frac{H-1}{N} I \right)+f\left(\frac{W-1}{N} I \right) \right] (N-[I]) .
$$
Thus equation \eqref{(5.14)} can be approximated as
\begin{equation}\label{csaladode45}
\dot{I}=\tau(N-I)\Bigg[f\Bigg(\frac{H-1}{N} I \Bigg)+f\Bigg(\frac{W-1}{N} I \Bigg)\Bigg]-\gamma I,
\end{equation}
The solution of this equation is compared to simulation in Figure~\ref{fig:csaladok} for two different hypergraphs. We can observe that the mean-field approximation gives better agreement when the hypergraph is homogeneous, i.e. for the case $H=10$, $W=10$.

\begin{figure}[h!]
\begin{center}
\includegraphics[scale=0.5]{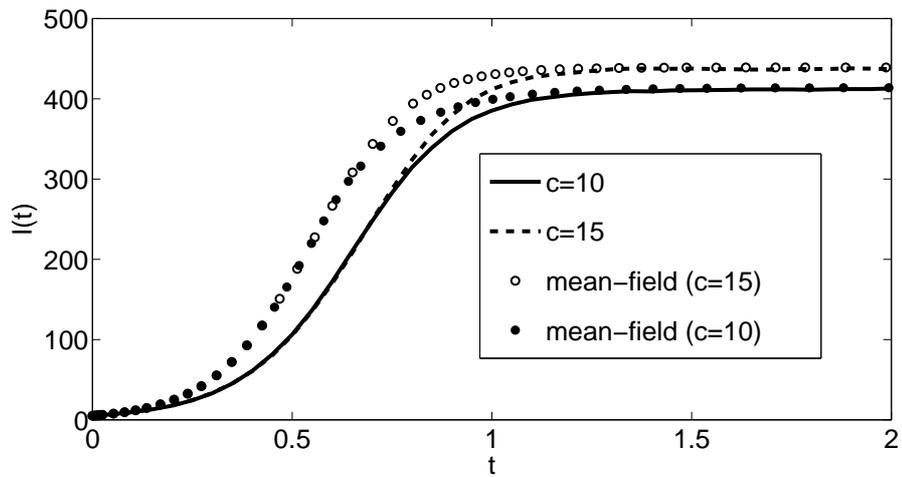}
\caption{Time dependence of the number of infected nodes for a regular random hypergraph, in which every node belongs to $d=16$ hyperedges, each of which is of size $e=20$. The simulation curves belonging to $c=10$ and $c=15$ in \eqref{eq:f} are shown together with the prevalence given by the mean-field equation \eqref{veletlenode45}. The parameter values are $N=500$, $M=400$, $\gamma =1$ and $\tau=0.03$. } \label{fig:meanf_veletlen2}
\end{center}
\end{figure}

Consider now regular random hypergraphs, in which each node belongs to $d$ hyperedges and each hyperedge is of size $e$. Denoting the number of infected nodes by $I$, the average number of infected neighbours of a node in a hyperedge can be approximated by $(e-1)\frac{I}{N} $. Thus the second summation in \eqref{NSIf} consists of a single term and can be approximated as $d f\left((e-1)\frac{I}{N}\right)$. Since this is independent of $l$, the double sum in \eqref{NSIf} reduces to
$$
N_{SI}^f(S_j^k) \approx (N-k) d f\left(\frac{e-1}{N} I \right),
$$
because the number of susceptible nodes in state $S_j^k$ is $N-k$ and each node belongs to $d$ hyperedges. Now \eqref{(5.12)} leads to the approximation
$$
[SI]\approx d f\left(\frac{e-1}{N} I \right) \sum_{k=0}^N \sum_{j=1}^{c_k} (N-k) X_j^k(t) =  d f\left(\frac{e-1}{N} I \right) (N-[I]) .
$$
Thus, for regular random hypergraphs, equation \eqref{(5.14)} can be approximated as
\begin{equation}\label{veletlenode45}
\dot{I}=\tau(N-I)\, d\, f\Bigg(\frac{e-1}{N} I \Bigg) -\gamma I.
\end{equation}
The solution of this equation is compared to simulation in Figure~\ref{fig:meanf_veletlen2} for two different values of $c$. We can see that for regular random hypergraphs the mean-field approximation performs well and for the steady state it gives excellent agreement. Similarly to the case of networks, the solution of the mean-field equation increases faster than the simulated curve, but their steady states are close to each other.

%%%%%%%%%%%%%%%%%%%%%%%%%%%%%%%%%%%%%%%%%%%%%%%%%%%%%%%%%%%%%%%%%%%%%%%%%%%%%%%%%%%%%%%%%%
\section{Discussion}
%%%%%%%%%%%%%%%%%%%%%%%%%%%%%%%%%%%%%%%%%%%%%%%%%%%%%%%%%%%%%%%%%%%%%%%%%%%%%%%%%%%%%%%%%%

In this paper the methods of mathematical modeling of epidemic propagation on networks is extended to hypergraphs. The aim of this extension is to account for both the community structure and the nonlinear dependence of the infection pressure on the number of infected neighbours. The novelty of the model is that a susceptible individual is assumed to become infected with probability $1-\exp(-r \Delta t)$ in a small time interval $\Delta t$ with rate $r=\tau \sum_{h} f(k_h)$, where the summation is for those hyperedges $h\in \mathcal{E}$ that contain the susceptible node, $k_h$ denotes the number of infected nodes in the hyperedge $h$ and $f$ is a function given in \eqref{eq:f}. The simulation algorithm, developed for networks, is extended to hypergraphs to account for this new transition probability function. Individual-based stochastic simulations were run on three different types of hypergraphs, configuration random hypergraphs, hypergraphs with hyperedges created from the cliques of a power-law random graph and bi-uniform hypergraphs the vertices of which are divided into households and workplaces randomly. The effects of hypergraph structure and the model parameters are investigated via individual-based simulation results. The exact master equations of the epidemic spreading are derived for an arbitrary hypergraph given by its incidence matrix. Based on these, moment closure approximation and mean-field models are introduced and compared to individual-based stochastic simulations.

The paper is the first step in extending the mathematical modeling of epidemic spreading from networks to hypergraphs. There are several directions where this extension can be continued. One of these is to develop and investigate pairwise models, in which not only the expected value of susceptible and infected nodes, but also the extension of pairs, as introduced in \eqref{(5.12)}, are determined. The exact master equation enables us to use lumping for reducing the size of the system when the network has a special symmetry. The idea of lumping is also extendable to the case of hypergraphs. Our results show that the simple mean-field models developed in Section \ref{sec:MF} performs well for homogeneous hypergraphs. The investigation of heterogeneous hypergraphs and the derivation and investigation of the corresponding heterogeneous mean-field models is also a challenging subject. The functional form and parameters of $f$ in \eqref{eq:f} could be determined based on real epidemic propagation data by using some fitting procedure. This approach is beyond the scope of this paper and may be the subject of future work.

%%%%%%%%%%%%%%%%%%%%%%%%%%%%%%%%%%%%%%%%%%%%%%%%%%%%%%%%%%%%%%%%%%%%%%%%%%%%%%%%%%%%%%%%%%
\section*{Acknowledgements} P\'eter L. Simon acknowledges support from Hungarian Scientific Research Fund, OTKA, (grant no. 115926).
Gyula Y. Katona acknowledges support from OTKA (grant no. 108947).
%%%%%%%%%%%%%%%%%%%%%%%%%%%%%%%%%%%%%%%%%%%%%%%%%%%%%%%%%%%%%%%%%%%%%%%%%%%%%%%%%%%%%%%%%%

%\end{document}
\newpage

\end{document}